\documentclass{article}
\usepackage{amsmath}
\usepackage{amsfonts,amssymb,latexsym,epsfig,amsthm}
\usepackage{epsfig}

\newtheorem{thm}{Theorem}

\DeclareMathOperator{\K}{\mathcal{K}}

\newcommand{\A}{\mathbf{A}}

\DeclareMathOperator{\Tr}{Tr}

\newcommand{\Prob}{\mathbb{P}}
\newcommand{\C}{\mathbb{C}}
\newcommand{\R}{\mathbb{R}}

\DeclareMathOperator{\diag}{diag}

\newcommand{\X}{\mathbf{X}}

\newcommand{\Z}{\mathbf{Z}}

\newcommand{\J}{\mathbf{J}}
\newcommand{\M}{\mathbf{M}}

\newcommand{\B}{\mathbf{B}}
\newcommand{\D}{\mathbf{D}}
\newcommand{\E}{\mathbf{E}}
\newcommand{\F}{\mathbf{F}}
\newcommand{\Con}{\mathbf{C}}
\newcommand{\U}{\mathbf{U}}

\newcommand{\Q}{\mathbf{Q}}
\newcommand{\SC}{\mathbf{S}}

\newcommand{\Cauchy}{\mathcal{C}^1}

\newcommand{\I}{\mathbb{I}}

\newcommand{\vecw}{\mathbf{w}}
\newcommand{\vecv}{\mathbf{v}}
\newcommand{\vect}{\mathbf{t}}
\newcommand{\vecx}{\mathbf{x}}

\newcommand{\e}{\epsilon}
\newcommand{\de}{\delta}

\usepackage{color}
\definecolor{blue}{rgb}{0,0,1}
\definecolor{red}{rgb}{1,0,0}

\title{\bf On the Christoffel-Darboux kernel for random Hermitian matrices with external source}

\author{\textbf{Jinho Baik}
\footnote{Department of Mathematics, University of Michigan, Ann
Arbor, MI, 48109, baik@umich.edu}}

\date{\today}

\begin{document}
\maketitle

\begin{abstract}
Bleher and Kuijlaars, and Daems and Kuijlaars showed that the correlation functions of 
the eigenvalues of a random matrix from unitary ensemble with external source can be expressed in terms of the Christoffel-Darboux kernel for multiple orthogonal polynomials. 
We obtain a representation of this Christoffel-Darboux kernel in terms of the usual orthogonal polynomials.  
\end{abstract}

\section{Introduction}

Fix a Hermitian matrix $A$, and a function $V(x)$ on $\R$, 
which decays sufficiently fast.
Consider the ensemble of $n\times n$ Hermitian matrices with the density function defined by 
\begin{equation}\label{eq:d1}
 p(M)=c_n\cdot e^{-\Tr(V(M)-AM)}.
\end{equation}
Here $c_n$ is the normalization constant,
\begin{equation}
    c_n= \bigg( \int e^{-\Tr(V(M)-AM)} dM \bigg)^{-1}, 
\end{equation}
where the integral is over all $n\times n$ Hermitian matrices. Due to the unitary-invariance of the density function, we can assume without loss of generality that $A$ is a diagonal matrix:
\begin{equation}
	A= \diag(a_1, a_2, \cdots, a_n)
\end{equation}
where $a_1\ge a_2\ge \cdots \ge a_n$. When $A=0$, the above ensemble is the standard unitary-invariant ensemble in the theory of random matrices \cite{Mehta}. The case when $A\neq 0$ is called a random matrix model with external source, and has been studied in, for example, 
\cite{BrezinHikami, ZinnJustin1, ZinnJustin2, BleherKuijlaars1, BleherKuijlaars2, BleherKuijlaars3}. 
It is well known that, using the famous Harish-Chandra-Itzykson-Zuber integral \cite{HC, IZ},~\eqref{eq:d1} induces the probability density function of the eigenvalues $\lambda_1\ge  \cdots \ge\lambda_n$ of $M$ given by 
\begin{equation}
   p(\lambda)=C_n \cdot
    \det(e^{a_i\lambda_j})_{1\le i,j\le n}
    \prod_{1\le i<j\le n} (\lambda_j-\lambda_i)
 \prod_{j=1}^n e^{-V(\lambda_j)}
\end{equation}
where $C_n=C_n(a_1, \cdots, a_n)$ is the (different) normalization constant.

In random matrix theory, one of the fundamental object of interest is the correlation functions of the eigenvalues. When $A=0$, the correlation functions are expressible in terms of the
Christoffel-Darboux kernel of the orthogonal polynomials with respect to 
the weight $e^{-V(x)}$ \cite{Mehta}. Using the Riemann-Hilbert problem (RHP) representation of orthogonal polynomials \cite{FIK}, the asymptotics of orthogonal polynomials with respect to a general class of weights have been obtained \cite{DKMVZ2, DKMVZ3}, hence yielding limit theorems for correlation functions.  

When $A\neq 0$, it was shown by Bleher and Kuijlaars \cite{BleherKuijlaars2004} that the correlation functions of the eigenvalues can be expressed in terms of a 
kernel involving a sum of multiple orthogonal polynomials.
This summation was further simplified by Daems and Kuijlaars 
\cite{DaemsKuijlaars}, yielding a Christoffel-Darboux formula for multiple orthogonal polynomials. 
The purpose of this paper is to re-express the Christoffel-Darboux kernel obtained by Daems and Kuijlaars in terms of orthogonal polynomials.

In order to state the result, we first discuss the result of \cite{DaemsKuijlaars}. 
Assume that $n-r$ eigenvalues of $A$ are zero. If $A$ has no zero eigenvalue, then we can simply shift $A$ by $A-a_1\I$ and $V(x)$ by $V(x)-a_1x$. Assume further that all non-zero eigenvalues of $A$ are simple: 
\begin{equation}\label{eq:A1}
	A=\diag(a_1, a_2, \cdots, a_r, 0, 0, \cdots, 0)
\end{equation}
where $a_1>a_2>\cdots >a_r$ and $a_j\neq 0$ for all $j=1, \cdots, r$. 
Where there are multiple eigenvalues, one obtain obtain the relevant result by taking limits such as $a_1\to a_2$ in the algebraic formulas below.

Let $\Z$ be the solution to the following $(r+2)\times (r+2)$ RHP:
\begin{itemize}
\item $\Z(z)$ is analytic for $z\in \C\setminus \R$ and is continuous up to the boudary,
\item $\Z_+(z)=\Z_-(z)\J_{\Z}(z)$ for $z\in \R$ where
\begin{equation}\label{eq:RHPforZ}
	\J_{\Z}(z)= \begin{pmatrix}
	1 & e^{-V(z)} & e^{-V(z)+a_1z} & \cdots &e^{-V(z)+a_rz} \\
	0 &1& 0 &\cdots & 0 \\
	0 & 0 & 1 &\cdots & 0\\
	\vdots &\vdots &  & \ddots &  \\
	0 & 0 & 0 &\cdots & 1\\
	\end{pmatrix},
\end{equation}

\item and\footnote{Throughout the paper $\I_k$ denotes the $k\times k$ identity matrix.}
\begin{equation}\label{eq:Zasymp}
\Z(z)= (\I_{r+2}+O(z^{-1})) \begin{pmatrix}
	z^n&0  &0&\cdots &0 \\
	0 & z^{-n+r}  &0&\cdots &0 \\
	0 & 0 & z^{-1}  &\cdots &0\\
	\vdots &\vdots &  & \ddots & \\
	0 & 0 & 0 &\cdots & z^{-1}\\
	\end{pmatrix}
\end{equation}
as $z\to \infty$.
\end{itemize}
Here for $z\in \R$, $\Z_{+}(z)$ (resp. $\Z_-(z)$) is the limit of $\Z(w)$ as $w\to z$ where $w$ is in the upper (resp. lower) half plane.
It is known that the solution $\Z$ exists, is unique, and is expressed in terms of multiple orthogonal polynomials of type II \cite{MultiOPRHP}. 
Define the kernel 
\begin{equation}\label{eq:kernel}
\begin{split}
	& \K(x,y) =\frac{e^{-V(x)}}{2\pi i(x-y)}
\begin{pmatrix}0 & 1 & e^{a_1y} & \cdots & e^{a_ry}\end{pmatrix}
\Z_{\pm}(y)^{-1}\Z_{\pm}(x) \begin{pmatrix}
1 \\ 0 \\ \vdots \\ 0
\end{pmatrix}.
\end{split}
\end{equation}
Note that since $(1,0,\cdots, 0)^t$ is a right eigenvector of $\J_{\Z}(x)$ corresponding to the eigenvalue $1$, and $(0,1, e^{a_1y}, \cdots, e^{a_ry})$ is  a left eigenvector of $\J_{\Z}(y)$ corresponding to the eigenvalue $1$, the kernel $\K(x,y)$ is independent of the choice of the plus or the minus limits of $\Z(x)$ and $\Z(y)$.
Daems and Kuijlaars \cite{DaemsKuijlaars} showed that the correlation functions of 
the eigenvalues of random Hermitian matrix with density~\eqref{eq:d1} when $A$ is of form~\eqref{eq:A1} 
is expressible in terms of $\K$. For example, the largest eigenvalue has the distribution 
\begin{equation}
	\Prob(\lambda_{\max}\le s)= \det(1-\K|_{(s,\infty)})
\end{equation}
where $\K|_{(s,\infty)}$ is the operator on $L^2((s,\infty))$ whose kernel $\K(x,y)$.

Let us compare the kernel $\K$ with the Christoffel-Darboux kernel for orthogonal polynomials. 
Let $\K_0$ be $\K$ when $A=0$, i.e. without any external source. 
Concretely, $\K_0$ is expressed in terms of the following RHP. Let $\X$ be the solution to the following $2\times 2$ RHP: 
\begin{itemize}
\item $\X(z)$ is analytic for $z\in \C\setminus \R$ and is continuous up to the boudary,
\item $\X_+(z)=\X_-(z)\J_{\X}(z)$ for $z\in \R$ where
\begin{equation}\label{eq:RHPforX}
	\J_{\X}(z)= \begin{pmatrix}
	1 & e^{-V(z)}  \\
	0 &1&
	\end{pmatrix},
\end{equation}

\item and
\begin{equation}
\X(z)= (\I_2+O(z^{-1})) \begin{pmatrix}
	z^n&0 \\
	0 & z^{-n}
	\end{pmatrix}
\end{equation}
as $z\to \infty$.
\end{itemize}
This is precisely the RHP for orthogonal polynomials \cite{FIK}. Let $\pi_n(z)=z^n+\cdots$ be the monic orthogonal polynomial of degree $n$ with respect to the weight $e^{-V(x)}$:
\begin{equation}
    \int_{\R} \pi_n(x) \pi_m(x) e^{-V(x)}dx = 0, \qquad n\neq m.
\end{equation}
Let $\kappa_n$ be the positive constant such that $\{\kappa_n\pi_n(z)\}$ becomes a sequence of orthonormal polynomials.
Then \cite{FIK},
\begin{equation}
    \X(z)= \begin{pmatrix}
    \pi_n(z) &  *  \\
    -2\pi i\kappa_{n-1}^2 \pi_{n-1}(z) &  *
    \end{pmatrix},
\end{equation}
where the second column is given by the Cauchy transform of the first column times the weight $e^{-V(x)}$. However we do not need the explicit form of the second column throughout this paper. 
Then\footnote{Usually $\K_0$ is defined by the term $e^{-V(x)}$ replaced by $e^{-\frac12(V(x)+V(y))}$. But they are conjugate of each other, and have the same spectrum.}  
\begin{equation}\label{eq:kernel0}
\begin{split}
	\K_0(x,y)
    &= \frac{e^{-V(x)}\kappa_{n-1}^2\big(\pi_n(x)\pi_{n-1}(y)-\pi_{n-1}(x)\pi_n(y)\big)}{x-y} \\
    &=\frac{e^{-V(x)}}{2\pi i(x-y)}
\begin{pmatrix}0 & 1 \end{pmatrix}
\X_+(y)^{-1}\X_+(x) \begin{pmatrix}
1 \\ 0
\end{pmatrix}.
\end{split}
\end{equation}

\bigskip

Now we state the result of this paper. 
Define the column vector of length $r$,
\begin{equation}\label{eq:vectdef}
\begin{split}
	\vect(z) = \begin{pmatrix}
    \pi_{n-r}(z) & \pi_{n-r+1}(z) & \cdots & \pi_{n-1}(z)
	\end{pmatrix}^t.
\end{split}
\end{equation}
Also define the column vector
\begin{equation}\label{eq:vdef}
\begin{split}
	\vecv(z)= \begin{pmatrix}
	e^{a_1z} & e^{a_2 z} &\cdots &e^{a_r z}	\end{pmatrix}^t
\end{split}
\end{equation}
and set
\begin{equation}\label{eq:wdef}
    \vecw(z)^t= \vecv(y)^t
	-  \int_{-\infty}^\infty \K_0(s,y) \vecv(s)^t ds.
\end{equation}
Define an $r\times r$ matrix
\begin{equation}\label{eq:defB}
\begin{split}
	\B= \int_{-\infty}^\infty \vect(s)\vecv(s)^{t} e^{-V(s)}ds.
\end{split}
\end{equation}

\begin{thm}\label{thm:main}
Let $A$ be as in~\eqref{eq:A1}. Then 
the matrix $\B$ is invertible and
\begin{equation}~\label{eq:main}
\begin{split}
	\big(\K(x,y)-\K_0(x,y)\big)e^{V(x)}  
= \vecw(y)^t\B^{-1} \vect(x).
\end{split}
\end{equation}
\end{thm}

\medskip

From a standard linear algebra computation (see~\eqref{eq:AintBigA} below), the above result can also be written as 
\begin{equation}
\begin{split}
	\big(\K(x,y)-\K_0(x,y)\big)e^{V(x)}  
	= -\frac{1}{\det \B}
	\det\begin{pmatrix}
	0 & \vecw(y)^t\\
	\vect(x) & \B
	\end{pmatrix}.
\end{split}
\end{equation}

\bigskip

Note that 
the RHP~\eqref{eq:RHPforZ} for multiple orthogonal polynomials are of size higher than $2$, contrary to the RHP~\eqref{eq:RHPforX} for orthogonal polynomials which is of size $2$. There have been many exciting recent progresses in solving $k\times k$ RHP asymptotically when $k=3,4$ (see for example, \cite{BleherKuijlaars1, BleherKuijlaars2, BleherKuijlaars3, DaemsKuijlaarsFour, DuitsKuijlaarsTwoMatrix}). Nevertheless the asymptotic study of RHP of general higher size is still a challenge.
The main motivation of this paper is to find an expression of the correlation function involving only $2\times 2$ RHP's. 

A special case is the so-called spiked model of random Hermitian matrices with external source when $n$ tends to infinity but $r$ is kept fixed.
In this case, the formula~\eqref{eq:main} seems well suited for the asymptotic analysis since the right-hand-side of~\eqref{eq:main} is of fixed rank $r$. In order to control this term asymptotically, one needs the first $r$ sub-leading terms of orthogonal polynomials in the asymptotic expansion explicitly. This requires a refinement of the analysis \cite{DKMVZ2, DKMVZ3}.
In the future work, we plan to work on this problem for general weight $V(x)$.
The particular case when $V(x)=x^2$ and also the Laguerre version of the problem have been studied in \cite{BBP, Peche} in the context of sample covariance matrices, which obtained various limit theorems and interesting phase transition phenomenon. The goal of the future work is to prove the universality for spiked model with general weight $V$.
On the other hand, when the rank $r$ of $A$ grows with $n$,~\eqref{eq:main} does not seem to be suitable for asymptotic analysis as the rank of the right-hand-side of~\eqref{eq:main} grows with $n$. For such case, one needs to analyze of the RHP~\eqref{eq:RHPforZ} for multiple orthogonal polynomials directly, and \cite{BleherKuijlaars1, BleherKuijlaars2, BleherKuijlaars3} worked out such asymptotic analysis when $V(x)=x^2$, and the half of the eigenvalues of $A$ are zero and the other half are the same, say $a$.

The proof of the theorem is based on algebraic manipulations of RHP's.

\medskip

\noindent {\it Acknowledgments: 
The work of the authors was supported in part by NSF grants DMS0457335, DMS075709. }

\section{Proof of Theorem~\ref{thm:main}}

Set
\begin{equation}\label{eq:Mdef}
\begin{split}
	\M(z):= \Z(z)\begin{pmatrix}
	\X(z)^{-1} & 0 \\
	0 & \I_r
	\end{pmatrix}.
\end{split}
\end{equation}
Then $\M(z)$ solves a new RHP. Let $\X_{ij}$ denote the $(i,j)$-entry of $\X$ and set
\begin{equation}
	\vecx(z)= \begin{pmatrix}
	\X_{11}(z) \\ \X_{21}(z)
	\end{pmatrix}
	= \begin{pmatrix}
	\pi_n(z) \\ -2\pi i\kappa_{n-1}^2\pi_{n-1}(z)
	\end{pmatrix}.
\end{equation}
The matrix $\M(z)$ satisfies:
\begin{itemize}
\item $\M(z)$ is analytic for $z\in \C\setminus \R$ and is continuous up to the boudary,
\item $\M_+(z)= \M_-(z)\J_{\M}(z)$ for $z\in \R$, where
\begin{equation}
\begin{split}
	\J_{\M}(z)
	&=\begin{pmatrix}
	1 & 0 & \X_{11}(z)e^{-V(z)+a_1z} & \cdots &\X_{11}(z)e^{-V(z)+a_rz} \\
	0 &1& \X_{21}(z)e^{-V(z)+a_1z} & \cdots &\X_{21}(z)e^{-V(z)+a_rz} \\
	0 & 0 & 1 &\cdots & 0\\
	\vdots &\vdots &  & \ddots &  \\
	0 & 0 & 0 &\cdots & 1\
	\end{pmatrix}   \\
	&= \begin{pmatrix}
	\I_2 & \vecx(z)\vecv(z)^t e^{-V(z)}\\
	0 & \I_r
	\end{pmatrix},
\end{split}
\end{equation}
\item
and as $z\to \infty$,
\begin{equation}\label{eq:Masz}
\begin{split}
	\M(z)
	&=\big(\I_{r+2}+O(z^{-1})\big)
	\begin{pmatrix}
	z^n\X_{22}(z) & -z^n\X_{12}(z) & 0 & \cdots &0 \\
	-z^{-n+r}\X_{21}(z) & z^{-n+r}\X_{11}(z)& 0 & \cdots &0 \\
	0 & 0 & z^{-1} &\cdots & 0\\
	\vdots &\vdots &  & \ddots &  \\
	0 & 0 & 0 &\cdots & z^{-1}
	\end{pmatrix}.
\end{split}
\end{equation}
\end{itemize}
Here $\I_{r+2}+O(z^{-1})$ in~\eqref{eq:Masz} and $\I_{r+2}+O(z^{-1})$ in~\eqref{eq:Zasymp} are the same function.
It will be clear in the below that the second column of this term plays a distinct role, and we denote it by 
\begin{equation}\label{eq:edeltadef}
    (0,1, 0, \cdots, 0)^t+ (\e_1(z), \e_2(z), \de_1(z), \cdots, \de_r(z))^t .
\end{equation}
Since the term $\I_{r+2}+O(z^{-1})$ equals $\Z(z) \diag( z^{n}, z^{-n+r}, z^{-1}, \cdots, z^{-1})^{-1}$ and $\Z(z)$ is given by the multiple orthogonal polynomials of type II \cite{MultiOPRHP},  the term $O(z^{-1})$ indeed have an asymptotic series in the powers of $z^{-1}$. Let
\begin{equation}\label{eq:ede}
\begin{split}
    \e_k(z)&= \e_k^1z^{-1} + \e_k^2 z^{-2} + \cdots + \e_k^r z^{-r} + O(z^{-r-1}), \qquad k=1,2,\\
    \de_k(z) &= \de_k^1z^{-1} + \de_k^2 z^{-2} + \cdots + \de_k^r z^{-r} + O(z^{-r-1}), \qquad k=1,\cdots, r,
\end{split}
\end{equation}
as $z\to\infty$ for some constants $\e_k^j$, $\de_\ell^j$, $j=1, \cdots, r$, $k=1,2$ and $\ell=1,\cdots, r$. 

We solve $\M$ in terms of $\X$. Write $\M$ in a block-form: let
\begin{equation}\label{eq:MblockJ}
\begin{split}
	\M(z)= \begin{pmatrix}
	\M^{(11)}(z) & \M^{(12)}(z) \\
	\M^{(21)}(z) & \M^{(22)}(z)
	\end{pmatrix}
\end{split}
\end{equation}
where $\M^{(11)}$ denotes the upper-left $2\times 2$ block of $\M$, and $\M^{(22)}$ denotes the lower-right $r\times r$ block of $\M$. The jump condition on $\M$ then becomes 
\begin{equation}\label{eq:Mjumpblock}
\begin{split}
	\begin{pmatrix}
	\M^{(11)}_+(z) & \M^{(12)}_+(z) \\
	\M^{(21)}_+(z) & \M^{(22)}_+(z)
	\end{pmatrix}
	= 
	\begin{pmatrix}
	\M^{(11)}_-(z) & \M^{(11)}_-(z)\vecx(z)\vecv(z)^t e^{-V(z)}+  \M^{(12)}_-(z) \\
	\M^{(21)}_-(z) & \M^{(21)}_-(z)\vecx(z)\vecv(z)^t e^{-V(z)}+ \M^{(22)}_-(z)
	\end{pmatrix}
\end{split}
\end{equation}
for $z\in\R$. 
Hence $\M^{(11)}(z)$ is continuous on $\R$ and therefore is entire.
As $z\to\infty$, the RHP for $\X$ implies that 
$z^n\X_{22}(z)=1+O(z^{-1})$ and $z^n\X_{12}(z)=O(z^{-1})$, and hence 
\begin{equation}\label{eq:M1100}
\begin{split}
	\M^{(11)}(z)
	&= \big(\I_{2}+O(z^{-1})\big) \begin{pmatrix}
	z^n\X_{22}(z) & -z^n\X_{12}(z) \\
	-z^{-n+r}\X_{21}(z) & z^{-n+r}\X_{11}(z)
	\end{pmatrix} \\
	&= \begin{pmatrix}
	1+O(z^{-1}) & \e_1(z) \\ O(z^{-1}) & 1+\e_2(z)
	\end{pmatrix} \begin{pmatrix}
	z^n\X_{22}(z) & -z^n\X_{12}(z) \\
	-z^{-n+r}\X_{21}(z) & z^{-n+r}\X_{11}(z)
	\end{pmatrix} \\
	&= \begin{pmatrix}
	1 &  0\\
	0 & 0
\end{pmatrix} + O(z^{-1})
+ \begin{pmatrix}
	\e_1(z) \\ 1+\e_2(z)
\end{pmatrix}\begin{pmatrix}
	-z^{-n+r}\X_{21}(z), & z^{-n+r}\X_{11}(z)
\end{pmatrix}.
\end{split}
\end{equation}
Let $\{b_j\}$ and $\{c_j\}$ be the coefficients of $\X_{11}$ and $\X_{21}$:
\begin{equation}
\begin{split}
	&\X_{11}(z)=\pi_n(z) = b_0z^n+ b_1z^{n-1}+ b_2z^{n-2}+ \cdots + b_n, \\
	&\X_{21}(z)=-2\pi i\kappa_{n-1}^2\pi_{n-1}(z) = c_0z^n+ c_1z^{n-1}+ c_2z^{n-2}+ \cdots + c_n.
\end{split}
\end{equation}
Note that
\begin{equation}\label{eq:bc0}
    b_0=1 \qquad c_0=0.
\end{equation}
From the Liouville's theorem, $\M^{(11)}(z)$ is a polynomial. This polynomial is obtained by expanding out the expansions of $\epsilon_k(z)$ and $\X_{ij}(z)$ in~\eqref{eq:M1100}, and collecting only the terms of non-negative powers of $z$. Explicitly, 
set
\begin{equation}
\begin{split}
	P_n^j(z) &:= b_0z^j+ b_1z^{j-1}+ b_2z^{j-2}+ \cdots+ b_j \\
	Q_n^j(z) &:= c_0z^j+c_1z^{j-1}+ c_2z^{j-2}+ \cdots + c_j.
\end{split}
\end{equation}
Then
\begin{equation}\label{eq:M11}
\begin{split}
	\M^{(11)}(z)
	= &\D(z)
	+ \begin{pmatrix} \e_1^1 & \e_1^2 & \cdots & \e_1^r \\ \e_2^1 & \e_2^2 & \cdots & \e_2^r
	 \end{pmatrix}\E(z)
\end{split}
\end{equation}
where
\begin{equation}
    \D(z):= \begin{pmatrix} 1 & 0 \\
	 -Q_n^r(z) &	
	 P_n^r(z)
	 \end{pmatrix}
\end{equation}
and
\begin{equation}\label{eq:E}
    \E(z):= \begin{pmatrix}
	 -Q_n^{r-1}(z) & P_n^{r-1}(z) \\
	 -Q_n^{r-2}(z) & P_n^{r-2}(z) \\
	 \vdots & \vdots  \\
	 -Q_n^{0}(z) & P_n^{0}(z)
	 \end{pmatrix}.
\end{equation}

The jump condition for $\M^{(12)}$ is (see~\eqref{eq:Mjumpblock})
\begin{equation}
\begin{split}
	\M^{(12)}_+(z) = \M^{(12)}_-(z) + \M^{(11)}(z) \vecx(z) \vecv(z)^te^{-V(z)}, \qquad z\in\R.
\end{split}
\end{equation}
Here $\M^{(11)}_-=\M^{(11)}$ since $\M^{(11)}(z)$ is entire.
Using the Plemelj's formula for the solution of an additive RHP, we obtain
\begin{equation}\label{eq:M12Pl}
\begin{split}
	\M^{(12)}(z)= \frac1{2\pi i}
	\int_{\R} \frac1{s-z}\M^{(11)}(s) \vecx(s) \vecv(s)^te^{-V(s)}ds,
\end{split}
\end{equation}
for $z\in \C\setminus\R$. Note that there is no additive term of an entire function in this formula due to the asymptotics
\begin{equation}
\begin{split}
	\M^{(12)}(z)=O(z^{-2}), \qquad z\to \infty, 
\end{split}
\end{equation}
which follows from~\eqref{eq:Masz}.
This asymptotic condition, applied to~\eqref{eq:M12Pl}, further implies that 
\begin{equation}
\begin{split}
	\int_{\R} \M^{(11)}(s) \vecx(s) \vecv(s)^te^{-V(s)}ds =0.
\end{split}
\end{equation}
From~\eqref{eq:M11}, the equations become
\begin{equation}\label{eq:eqfore}
\begin{split}
    &\begin{pmatrix} \e_1^1 & \e_1^2 & \cdots & \e_1^r \\ \e_2^1 & \e_2^2 & \cdots & \e_2^r
	 \end{pmatrix} \Q
    = \Con,
\end{split}
\end{equation}
where
\begin{equation}\label{eq:Q}
    \Q:= \int_{\R} \E(s)\vecx(s)\vecv(s)^t e^{-V(s)}ds.
\end{equation}
and
\begin{equation}
    \Con:= -\int_{\R}
    \D(s)
	 \vecx(s)\vecv(s)^t e^{-V(s)}ds.
\end{equation}
We will see~\eqref{eq:eqfordel} below that $\Q$ is invertible. Hence~\eqref{eq:eqfore2} and~\eqref{eq:M11} imply that 
\begin{equation}\label{eq:eqfore2}
\begin{split}
    \M^{(11)}(z)=\D(z)+ \Con\Q^{-1}\E(z). 
\end{split}
\end{equation}

Now consider $\M^{(21)}$. The jump condition~\eqref{eq:Mjumpblock} implies that $\M^{(21)}(z)$ is an entire function. As $z\to \infty$, from~\eqref{eq:Masz}, 
we have
\begin{equation}
\begin{split}
	\M^{(21)}(z)
	&= \begin{pmatrix}
	O(z^{-1}) & \delta_1(z) \\
	O(z^{-1}) & \delta_2(z) \\
	\vdots & \vdots \\
	O(z^{-1}) & \delta_r(z) \\
\end{pmatrix}
\begin{pmatrix}
	z^n\X_{22}(z) & -z^n\X_{12}(z) \\
	-z^{-n+r}\X_{21}(z) & z^{-n+r}\X_{11}(z)
\end{pmatrix} \\
	&= O(z^{-1})+
\begin{pmatrix}
	\de_1(z) \\ \de_2(z) \\ \vdots \\ \de_r(z)
\end{pmatrix}\begin{pmatrix}
	-z^{-n+r}\X_{21}(z), & z^{-n+r}\X_{11}(z)
\end{pmatrix}.
\end{split}
\end{equation}
From the Liouville's theorem, $\M^{(21)}(z)$ is a polynomial of degree at most $r-1$, and is given by (recall~\eqref{eq:ede})
\begin{equation}\label{eq:M21}
\begin{split}
	\M^{(21)}(z)
	= \begin{pmatrix} \de_1^1 & \de_1^2 & \cdots & \de_1^r \\ \de_2^1 & \de_2^2 & \cdots & \de_2^r \\
    \vdots & \vdots & \ddots &\vdots \\ \de_r^1 & \de_r^2 & \cdots & \de_r^r
	 \end{pmatrix}
    \E(z),
\end{split}
\end{equation}
where $\E(z)$ is defined in~\eqref{eq:E}.

The jump condition for $\M^{(22)}$ is
\begin{equation}\label{eq:M22asymp}
\begin{split}
	\M^{(22)}_+(z) = \M^{(22)}_-(z) + \M^{(21)}(z) \vecx(z) \vecv(z)^te^{-V(z)}, \qquad z\in\R.
\end{split}
\end{equation}
Using the asymptotic condition
\begin{equation}
\begin{split}
	\M^{(22)}(z)=z^{-1}\I_r+ O(z^{-2}), \qquad z\to \infty,
\end{split}
\end{equation}
the Plemelj's formula implies that
\begin{equation}\label{eq:M22Pl}
\begin{split}
	\M^{(22)}(z)= \frac1{2\pi i}
	\int_{\R} \frac1{s-z}\M^{(21)}(s) \vecx(s) \vecv(s)^te^{-V(s)}ds, \quad z\in \C\setminus\R.
\end{split}
\end{equation}
The asymptotic condition~\eqref{eq:M22asymp} further implies that
\begin{equation}
\begin{split}
	\int_{\R} \M^{(21)}(s) \vecx(s) \vecv(s)^te^{-V(s)}ds = - 2\pi i \I_r.
\end{split}
\end{equation}
From~\eqref{eq:M21}, this implies that 
\begin{equation}\label{eq:eqfordel}
\begin{split}
    \begin{pmatrix} \de_1^1 & \de_1^2 & \cdots & \de_1^r \\ \de_2^1 & \de_2^2 & \cdots & \de_2^r \\
    \vdots & \vdots & \ddots &\vdots \\ \de_r^1 & \de_r^2 & \cdots & \de_r^r
	 \end{pmatrix} \Q
    =  -2\pi i\I_r
\end{split}
\end{equation}
where $\Q$ is in~\eqref{eq:Q}. This implies that $\Q$ is invertible, and 
\begin{equation}\label{eq:eqfordel2}
\begin{split}
    \M^{(21)}(z)= -2\pi i \Q^{-1}\E(z). 
\end{split}
\end{equation}

Therefore, combining~~\eqref{eq:M12Pl},~\eqref{eq:eqfore2},~\eqref{eq:M22Pl} and~\eqref{eq:eqfordel2}, the matrix $\M(z)$ is given by 
given by 
\begin{equation}\label{eq:Msol}
\begin{split}
    \M(z)
    &= \begin{pmatrix}
    \D(z) + \Con \Q^{-1} \E(z) & (\Cauchy \D)(z) + \Con\Q^{-1} (\Cauchy \E)(z) \\
    -2\pi i \Q^{-1} \E(z) & -2\pi i \Q^{-1} (\Cauchy\E)(z)
    \end{pmatrix} \\
    &= \begin{pmatrix} \I_2 & \Con \\ 0 & -2\pi i\I_r \end{pmatrix}
    \begin{pmatrix} \I_2 & 0 \\ 0 & \Q^{-1} \end{pmatrix}
    \begin{pmatrix} \D(z) & (\Cauchy \D)(z) \\ \E(z) & (\Cauchy \E)(z) \end{pmatrix}
\end{split}
\end{equation}
where we use the notation 
\begin{equation}
    (\Cauchy F)(z):= \frac1{2\pi i} \int_{\R} \frac1{s-z} F(s)\vecx(s) \vecv(s)^te^{-V(s)} ds, \qquad z\in \C\setminus \R
\end{equation}
for $2\times 2$ matrix-valued function $F(z)$ on $\R$.

Now consider the kernel~\eqref{eq:kernel}. From~\eqref{eq:Mdef},
\begin{equation}
\begin{split}
    &2\pi i\K(x,y)e^{V(x)}(x-y) \\
    &=\begin{pmatrix}0 & 1 & e^{a_1y} & \cdots & e^{a_ry}\end{pmatrix}
    \begin{pmatrix} \X_+(y)^{-1} & 0 \\ 0 & \I_r \end{pmatrix}
\M_+(y)^{-1}\M_+(x)
    \begin{pmatrix} \X_+(x) & 0 \\ 0 & \I_r \end{pmatrix}\begin{pmatrix}
1 \\ 0 \\ \vdots \\ 0
\end{pmatrix}  \\
    &= \begin{pmatrix} -\X_{21}(y) & \X_{11}(y) & e^{a_1y} & \cdots & e^{a_ry}\end{pmatrix}
\M_+(y)^{-1}\M_+(x)
    \begin{pmatrix} \X_{11}(x) \\ \X_{21}(x) \\ 0 \\ \vdots  \\0  \end{pmatrix}.
\end{split}
\end{equation}
On the other hand,  from~\eqref{eq:kernel0},
\begin{equation}
\begin{split}
    &2\pi i \K_0(x,y)e^{V(x)}(x-y) \\
    &= \begin{pmatrix} -\X_{21}(y) & \X_{11}(y) \end{pmatrix}
    \begin{pmatrix} \X_{11}(x) \\ \X_{21}(x) \end{pmatrix} \\
    &= \begin{pmatrix} -\X_{21}(y) & \X_{11}(y) & e^{a_1y} & \cdots & e^{a_ry}\end{pmatrix}
    \begin{pmatrix} \X_{11}(x) \\ \X_{21}(x) \\ 0 \\ \vdots  \\0  \end{pmatrix}.
\end{split}
\end{equation}
Hence
\begin{equation}
\begin{split}
    &2\pi i\big( \K(x,y)-\K_0(x,y)\big) e^{V(x)}\\
    &= \frac1{x-y}\begin{pmatrix} -\X_{21}(y) & \X_{11}(y) & e^{a_1y} & \cdots & e^{a_ry}\end{pmatrix}
    \big( \M_+(y)^{-1} \M_+(x) - \I_{r+2}\big)
    \begin{pmatrix} \X_{11}(x) \\ \X_{21}(x) \\ 0 \\ \vdots  \\0  \end{pmatrix} \\
    &= \begin{pmatrix} -\X_{21}(y) & \X_{11}(y) & e^{a_1y} & \cdots & e^{a_ry}\end{pmatrix}
    \M_+(y)^{-1}  \frac{\M(x) - \M(y)}{x-y}
    \begin{pmatrix} \X_{11}(x) \\ \X_{21}(x) \\ 0 \\ \vdots  \\0  \end{pmatrix}.
\end{split}
\end{equation}
Here in the last line, we replaced $\frac{\M_+(x) - \M_+(y)}{x-y}$ by $\frac{\M(x) - \M(y)}{x-y}$ 
since only the first two columns are needed when multiplied by the last vector, and the first two columns of $\M(z)$ are entire. 
Using the formula~\eqref{eq:Msol} for $\M$, we obtain
\begin{equation}
\begin{split}
    &2\pi i\big( \K(x,y)-\K_0(x,y)\big) e^{V(x)}\\
    &= \begin{pmatrix} -\X_{21}(y) & \X_{11}(y) & e^{a_1y} & \cdots & e^{a_ry}\end{pmatrix}
    \begin{pmatrix} \D(y) & (\Cauchy_+\D)(y) \\ \E(y) & (\Cauchy_+ \E)(y) \end{pmatrix}^{-1}
    \begin{pmatrix} \frac{\D(x)-\D(y)}{x-y} \\ \frac{\E(x)-\E(y)}{x-y} \end{pmatrix} \vecx(x).
\end{split}
\end{equation}
Here note that since $\M(z)$ is invertible and $\det \M(z)=1$, the last matrix in~\eqref{eq:Msol} is invertible and 
\begin{equation}
\begin{split}
    \det
    \begin{pmatrix} \D(z) & (\Cauchy_+\D)(z) \\ \E(z) & (\Cauchy_+ \E)(z) \end{pmatrix}
    =(-2\pi i)^{-r}\det \Q.
\end{split}
\end{equation}

Using the general identity
\begin{equation}\label{eq:AintBigA}
\begin{split}
   -\vecw_1 ^t \A^{-1}\vecw_2= \frac1{\det \A} \det \begin{pmatrix}
   0 &\vecw_1^t \\ \vecw_2 & \A
   \end{pmatrix}
\end{split}
\end{equation}
for an invertible matrix $\A$ and vectors $\vecw_1, \vecw_2$, we obtain
\begin{equation}
\begin{split}
    &2\pi i\big( \K(x,y)-\K_0(x,y)\big) e^{V(x)}\\
    &= -\frac{(-2\pi i)^r}{\det\Q}\det \begin{pmatrix} 0 &  \begin{pmatrix}-\X_{21}(y) & \X_{11}(y) & e^{a_1y} & \cdots & e^{a_ry} \end{pmatrix} \\
    \begin{pmatrix} \frac{\D(x)-\D(y)}{x-y} \\ \frac{\E(x)-\E(y)}{x-y} \end{pmatrix} \vecx(x)
     & \begin{pmatrix} \D(y) & (\Cauchy_+\D)(y) \\ \E(y) & (\Cauchy_+ \E)(y) \end{pmatrix} \end{pmatrix} \\
    &= -\frac{(-2\pi i)^r}{\det\Q}\det \begin{pmatrix} 0 &  \hat\vecx(y)  & \vecv(y)^t \\
    \frac{\D(x)-\D(y)}{x-y}\vecx(x)
     &  \D(y) & (\Cauchy_+\D)(y) \\ \frac{\E(x)-\E(y)}{x-y}\vecx(x) & \E(y) & (\Cauchy_+ \E)(y)  \end{pmatrix} 
\end{split}
\end{equation}
where
\begin{equation}
    \hat\vecx(y):= \begin{pmatrix} 0 & -1 \\ 1 &0 \end{pmatrix} \vecx(y)
    = \begin{pmatrix} -\X_{21}(y) \\ \X_{11}(y) \end{pmatrix}
\end{equation}
and 
\begin{equation}
\begin{split}
	\vecv(z)= \begin{pmatrix}
	e^{a_1z} & e^{a_2 z} &\cdots &e^{a_r z}	\end{pmatrix}^t
\end{split}
\end{equation}
as defined in~\eqref{eq:vdef}.
Set 
\begin{equation}
\begin{split}
    \Delta
    &:= \det \begin{pmatrix} 0 &  \hat\vecx(y)  & \vecv(y)^t \\
    \frac{\D(x)-\D(y)}{x-y}\vecx(x)
     &  \D(y) & (\Cauchy_+\D)(y) \\ \frac{\E(x)-\E(y)}{x-y}\vecx(x) & \E(y) & (\Cauchy_+ \E)(y)  \end{pmatrix}
\end{split}
\end{equation}
so that 
\begin{equation}\label{eq:KK0diff}
	2\pi i\big( \K(x,y)-\K_0(x,y)\big) e^{V(x)}= -\frac{(-2\pi i)^r}{\det \Q}\Delta.
\end{equation} 
We have 
\begin{equation}\label{eq:Delta3}
\begin{split}
    \Delta    
    &= \det \begin{pmatrix} 0 &  \hat\vecx(y)  & \vecv(y)^t \\
     \frac{\D(x)-\D(y)}{x-y}\vecx(x)
     &  \D(y) & (\Cauchy_+\D)(y) \\ \frac{\E(x)-\E(y)}{x-y}\vecx(x) & \E(y) & (\Cauchy_+ \E)(y)  \end{pmatrix}
     \begin{pmatrix} \I_1 & 0 & 0\\ 0 & \I_2 & - (\Cauchy_+ \I_2)(y) \\ 0 & 0& \I_r
     \end{pmatrix}\\
    &= \det \begin{pmatrix} 0 &  \hat\vecx(y)  & \vecv(y)^t -\hat\vecx(y)(\Cauchy_+ \I_2)(y)\\
    \frac{\D(x)-\D(y)}{x-y}\vecx(x)
     &  \D(y) & (\Cauchy_+\D)(y) -\D(y)(\Cauchy _+\I_2)(y) \\ \frac{\E(x)-\E(y)}{x-y}\vecx(x) & \E(y) & (\Cauchy_+ \E)(y)  -\E(y)(\Cauchy_+ \I_2)(y) \end{pmatrix} .
\end{split}
\end{equation}
As $\K_0(s,y)=\frac{e^{-V(s)}(\X_{11}(s)\X_{21}(y)-\X_{21}(s)\X_{11}(y))}{2\pi i(s-y)}=\frac{e^{-V(s)}\hat\vecx(y)^t\vecx(s)}{2\pi i(s-y)}$, the term in the $(13)$-block is
\begin{equation}
\begin{split}
    \vecv(y)^t -\hat\vecx(y)(\Cauchy_+ \I_2)(y)
    &= \vecv(y)^t - \frac{\hat\vecx(y)^t}{2\pi i} \int_{\R} \frac1{(s-y)_+}  \vecx(s) \vecv(s)^t e^{-V(s)}ds\\
    &= \vecv(y)^t - \int_{\R} \K_0(s,y) \vecv(s)^t ds = \vecw(y)^t
\end{split}
\end{equation}
(see~\eqref{eq:wdef}).
On the other hand, since the first row of $\D(z)$ is $(1,0)$ and
\begin{equation}
\begin{split}
    (\Cauchy_+\D)(y) -\D(y)(\Cauchy _+\I_2)(y) = \frac1{2\pi i} \int_{\R}
    \frac{\D(s)-\D(y)}{(s-y)_+} \vecx(s)\vecv(s)^t e^{-V(s)}ds,
\end{split}
\end{equation}
the second row of the matrix in~\eqref{eq:Delta3} equals $(0,1,0,0,\cdots, 0)$.
Consider the last row of the matrix. Since the last row of $\E(z)$ is $(-Q_n^0(z), P_n^0(z))=(-c_0, b_0)=(0,1)$, we find that the last row of the matrix in~\eqref{eq:Delta3} is $(0,0,1,0,\cdots, 0)$.
Hence the determinant of the matrix in~\eqref{eq:Delta3} is unchanged, up to the factor $(-1)^{r}$, if we remove the second and the last row (which is the $(r+3)$rd row), and the second and the third column. Therefore, noting that the second row of $\D(z)$ is same as the first row of $\E(z)$ with $r-1$ replaced by $r$, we arrive at the formula
\begin{equation}
\begin{split}
    \Delta= (-1)^{r}
    \det \begin{pmatrix} 0 &  \vecw(y)^t\\
    \frac{\F(x)-\F(y)}{x-y}\vecx(x)
    & (\Cauchy_+\F)(y) -\F(y)(\Cauchy _+\I_2)(y)
     \end{pmatrix}
\end{split}
\end{equation}
where
\begin{equation}
\begin{split}
	\F(z)
    &:=
	 \begin{pmatrix}
	 -Q_n^{r}(z) & P_n^{r}(z) \\
	 -Q_n^{r-1}(z) & P_n^{r-1}(z) \\
	 \vdots & \vdots  \\
	 -Q_n^{1}(z) & P_n^{1}(z)
	 \end{pmatrix}\\
    &= \begin{pmatrix} 1 & z &  \cdots &z^{r-1} & z^r \\
    0 & 1 & \cdots & z^{r-2} & z^{r-1} \\
    \vdots & \vdots & \ddots & \vdots &\vdots \\
     0 & 0 &\cdots & 1 & z
    \end{pmatrix}
	 \begin{pmatrix}
	 -c_r & b_r \\
	 -c_{r-1} & b_{r-1} \\
	 \vdots & \vdots  \\
	 -c_0 & b_0
	 \end{pmatrix}.
\end{split}
\end{equation}
The last line is an $r\times (r+1)$ matrix times an $(r+1)\times 2$ matrix.
Now
\begin{equation}
\begin{split}
    &\frac{\F(s)-\F(y)}{s-y}  \\
    &= \begin{pmatrix} 0&  1 & s+y &  \cdots &s^{r-1}+s^{r-2}y+ \cdots + y^{r-1} \\
    0& 0 & 1 & \cdots &  s^{r-2}+s^{r-3}y+ \cdots + y^{r-2}\\
    \vdots & \vdots & \vdots & \ddots & \vdots  \\
     0& 0 & 0 &\cdots & 1
    \end{pmatrix}
	 \begin{pmatrix}
	 -c_r & b_r \\
	 -c_{r-1} & b_{r-1} \\
	 -c_{r-2} & b_{r-2} \\
	 \vdots & \vdots  \\
	 -c_0 & b_0
	 \end{pmatrix} \\
    &= \begin{pmatrix} 1 & s+y &  \cdots &s^{r-1}+s^{r-2}y+ \cdots + y^{r-1} \\
    0 & 1 & \cdots &  s^{r-2}+s^{r-3}y+ \cdots + y^{r-2}\\
    \vdots & \vdots & \ddots & \vdots  \\
     0 & 0 &\cdots & 1
    \end{pmatrix}
	 \begin{pmatrix}
	 -c_{r-1} & b_{r-1} \\
	 -c_{r-2} & b_{r-2} \\
	 \vdots & \vdots  \\
	 -c_0 & b_0
	 \end{pmatrix} \\
    &= \U(y)\U(s)\begin{pmatrix}
	 -c_{r-1} & b_{r-1} \\
	 -c_{r-2} & b_{r-2} \\
	 \vdots & \vdots  \\
	 -c_0 & b_0
	 \end{pmatrix}
    = \U(y) \E(s)
\end{split}
\end{equation}
where
\begin{equation}
    \U(z):=
	 \begin{pmatrix} 1 & z &  \cdots & z^{r-1}  \\
    0 & 1 & \cdots & z^{r-2}  \\
    \vdots & \vdots & \ddots & \vdots \\
     0 & 0 &\cdots & 1
    \end{pmatrix}.
\end{equation}
Therefore we obtain
\begin{equation}
\begin{split}
    \Delta
    &= (-1)^{r} \det \begin{pmatrix}
    0 & \vecw(y)^t \\ \U(y)\E(x) \vecx(x)&
    \frac{\U(y) }{2\pi i} \int_{\R} \E(s) \vecx(s)\vecv(s)^t e^{-V(s)}ds
    \end{pmatrix} \\
    &= (-1)^{r} \det \begin{pmatrix}
    0 & \vecw(y)^t \\ \E(x) \vecx(x) &
    \frac1{2\pi i} \int_{\R} \E(s) \vecx(s) \vecv(s)^t e^{-V(s)}ds
    \end{pmatrix}\\
    &= (-1)^{r}  \det \begin{pmatrix}
    0 & \vecw(y)^t \\ \E(x) \vecx(x) &
    \frac1{2\pi i} \Q
    \end{pmatrix}
\end{split}
\end{equation}
since $\det \U(y)=1$.

Now we compute 
\begin{equation}
\begin{split}
    \E(z)\vecx(z)
    = \begin{pmatrix}
	-Q_n^{r-1}(z)\X_{11}(z)+P_n^{r-1}(z)\X_{21}(z)\\
	-Q_n^{r-2}(z)\X_{11}(z)+P_n^{r-2}(z)\X_{21}(z)\\
	\vdots \\
	-Q_n^{1}(z)\X_{11}(z)+P_n^{1}(z)\X_{21}(z)\\
	-Q_n^{0}(z)\X_{11}(z)+P_n^{0}(z)\X_{21}(z)
\end{pmatrix}.
\end{split}
\end{equation}
The polynomial
\begin{equation}
\begin{split}
    &-Q_n^j(z)\X_{11}(z)+P_n^j(z)\X_{21}(z)\\
    &= -(c_0z^j+\cdots + c_j)(b_0z^n+ \cdots + b_n)
    +(b_0z^j+ \cdots + b_j)(c_0z^n+ \cdots + c_n)\\
    &= -(c_0b_{j+1}-b_0c_{j+1})z^{n-1}+ \cdots
\end{split}
\end{equation}
is of degree $n-1$. Hence it is a linear combination of the orthogonal polynomials $\pi_0(z), \cdots, \pi_{n-1}(z)$. Define the real inner product, 
\begin{equation}
    <f,g> := \int_{\R} f(x)g(x) e^{-V(x)}dx.
\end{equation}
Recalling that $\X_{11}(z)=\pi_n(z)$ and $\X_{21}(z)=-2\pi i\kappa_{n-1}^2\pi_{n-1}(z)$, 
we have for $k=0,1,\cdots, n-j-2$,
\begin{equation}
    <-Q_n^j\X_{11}+P_n^j\X_{21}, \pi_k>
    = -<\pi_n, Q_n^j\pi_k> - 2\pi i\kappa_{n-1}^2 <\pi_{n-1}, P_n^j \pi_k>
    =0
\end{equation}
since $Q_n^j\pi_k$ is of degree $k+j-1$ (recall that $c_0=0$) and $P_n^j\pi_k$ is of degree $k+j$.
Also since $b_0=1$,
\begin{equation}
\begin{split}
    <-Q_n^j\X_{11}+P_n^j\X_{21}, \pi_{n-j-1}>
    &= -2\pi i\kappa_{n-1}^2<\pi_{n-1}, P_n^j\pi_{n-j-1}> \\
    &= -2\pi i\kappa_{n-1}^2<\pi_{n-1}, z^{n-1}>\\
    &= -2\pi i.
\end{split}
\end{equation}
Noting that $<\pi_k,\pi_k>=\kappa_k^{-2}$, we obtain
\begin{equation}\label{eq:Ex}
\begin{split}
    \E(z)\vecx(z)
    &= - 2\pi i \SC \vect(z)
\end{split}
\end{equation}
where
\begin{equation}
\begin{split}
   \SC=  \begin{pmatrix}
	\kappa_{n-r}^2 & * & \cdots & * & *\\
    0& \kappa_{n-r+1}^2 &\cdots & * & * \\
    \vdots & \vdots & \ddots & \vdots & \vdots \\
    0& 0& \cdots & \kappa_{n-2}^2 & *\\
    0& 0& \cdots & 0& \kappa_{n-1}^2
	\end{pmatrix}
\end{split}
\end{equation}
and $\vect(z)$ is defined in~\eqref{eq:vectdef}.
From~\eqref{eq:Q}, this implies that  
\begin{equation}\label{eq:QbyB}
    \Q= -2\pi i  \SC
    \B
\end{equation}
where $\B$ is defined in~\eqref{eq:defB}. 
As $\Q$ is invertible, so is $\B$. 
Thus
\begin{equation}
\begin{split}
	\Delta &= (-1)^{r} \det \begin{pmatrix} 0 & \vecw(y)^t \\
	-2\pi i \SC \vect(x) & -\SC \B
	\end{pmatrix} \\
	&= \det \SC \det \begin{pmatrix} 0 & \vecw(y)^t \\
	2\pi i  \vect(x) & \B
	\end{pmatrix} 
	= -2\pi i \det(\SC\B) \vecw(y)^t \B^{-1} \vect(x)
\end{split}
\end{equation}
using~\eqref{eq:AintBigA}. As $\det\Q=(-2\pi i)^r\det(\SC\B)$,  
we obtain from~\eqref{eq:KK0diff} that 
\begin{equation}
\begin{split}
	2\pi i \big( \K(x,y)-\K_0(x,y)\big) e^{V(x)}= 2\pi i\vecw(y)^t \B^{-1} \vect(x).
\end{split}
\end{equation}
This completes the proof of the Theorem. 

\bibliographystyle{plain}

\end{document}